\documentclass[11pt]{article}

\usepackage[T1]{fontenc}
\usepackage[utf8]{inputenc}
\usepackage[english]{babel}
\usepackage{algorithm}
\usepackage{algpseudocode}
\usepackage{doi}

\usepackage{tikz}
\usetikzlibrary{arrows.meta, positioning, shapes.geometric}

\usepackage[a4paper,margin=30mm]{geometry}
\usepackage{lmodern}
\usepackage{microtype}
\usepackage{setspace}
\usepackage{amsmath,amssymb}

\usepackage{booktabs}
\usepackage{tabularx}

\usepackage[numbers]{natbib}
\usepackage{placeins} % in preamble

\title{
	Dynamic Inclusion and Bounded Multi-Factor Tilts\\
	for Robust Portfolio Construction
}

\author{
	Roberto Garrone\\
	\vspace{0.5em}
	\small Preprint
}

\date{\today}

\begin{document}
	\maketitle
	
\begin{abstract}
This paper proposes a portfolio construction framework designed to remain robust under estimation error, non-stationarity, and realistic trading constraints. The methodology combines dynamic asset eligibility, deterministic rebalancing, and bounded multi-factor tilts applied to an equal-weight baseline.

Asset eligibility is formalized as a state-dependent constraint on portfolio construction, allowing factor exposure to adjust endogenously in response to observable market conditions such as liquidity, volatility, and cross-sectional breadth. Rather than estimating expected returns or covariances, the framework relies on cross-sectional rankings and hard structural bounds to control concentration, turnover, and fragility.

The resulting approach is fully algorithmic, transparent, and directly implementable. It provides a robustness-oriented alternative to parametric optimization and unconstrained multi-factor models, particularly suited for long-horizon allocations where stability and operational feasibility are primary objectives.
\end{abstract}

	\section{Introduction}
	
Mean--variance portfolio optimization is theoretically optimal under strong assumptions
regarding the stationarity and accurate estimation of expected returns and covariances.
In practice, these assumptions are routinely violated, leading to unstable portfolios,
excessive turnover, and poor out-of-sample performance \cite{michaud1989,demiguel2009,jagannathan2003}.

This paper presents a portfolio construction methodology that explicitly avoids these failure
modes. The proposed approach replaces parametric optimization with a rule-based allocation
that combines dynamic asset inclusion, an equal-weight baseline, and bounded multi-factor
tilts under a fixed rebalancing schedule. Rather than pursuing theoretical efficiency, the
objective is robustness, interpretability, and operational feasibility under estimation error,
non-stationarity, and realistic trading constraints. Robust alternatives to parametric optimization have been shown to deliver
more stable out-of-sample performance under realistic estimation error and
structural instability, particularly when portfolio constraints are made
explicit \cite{fabozzi2007, boyd2004}.

This framework should not be interpreted as a forecasting model, a regime-switching system, or an alpha-maximization strategy. It does not estimate time-varying expected returns, perform volatility targeting, or condition allocations on inferred market regimes. Adaptivity arises exclusively through deterministic eligibility constraints based on lagged and observable information, rather than through predictive signals or parameter switching.

The contribution of this work is the formalization of dynamic asset eligibility as a
state-dependent constraint on portfolio construction. By treating the investable universe
as endogenous to observable market conditions—such as liquidity, volatility, and
cross-sectional breadth—the framework allows factor exposure to adjust adaptively without
relying on parameter switching, regime detection, or forecast-driven optimization. The
resulting methodology is fully algorithmic and directly implementable, providing a robust
alternative to parametric optimization and unconstrained multi-factor models.

It is important to clarify what the proposed framework does not attempt to do. The method is
not designed to estimate time-varying expected returns, to perform tactical timing, or to
optimize portfolios with respect to volatility, tracking error, or benchmark-relative
objectives. Instead, it focuses on constructing a stable and interpretable allocation that
degrades gracefully under adverse conditions and adapts endogenously through eligibility
constraints rather than through explicit prediction or regime classification.

	\section{Design Principles}
	
	The proposed framework is guided by four design principles:
	
	\begin{enumerate}
		\item \textbf{Robustness over optimality}: preference is given to constructions that degrade gracefully under noise.
		\item \textbf{Cross-sectional dominance}: relative rankings are favored over level estimates\cite{grinold2000,qian2007,harvey2016}. 
		\item \textbf{Explicit constraints}: bounds and eligibility rules are structural, not emergent.
		\item \textbf{Operational realism}: rebalancing frequency, liquidity, and costs are integral inputs.
	\end{enumerate}
	
Cross-sectional ranking approaches are less sensitive to level estimation
	error and provide increased stability under non-stationarity, particularly
	when conditioning information is noisy or regime-dependent
	\cite{ferson1999, pedersen2015}. Taken together, these principles define a state-conditioned portfolio construction framework in which adaptivity emerges from endogenous universe selection rather than from forecast-driven optimization.

\section{Methodological Validity and Threats to Inference}
\label{sec:validity}

This section addresses methodological concerns commonly raised in empirical asset pricing
and portfolio construction studies, including endogeneity, look-ahead bias, data-snooping,
and factor redundancy. The objective is not to claim statistical optimality, but to demonstrate
that the proposed framework is internally coherent, information-consistent, and free from
structural biases that would invalidate empirical interpretation.

\subsection{Endogeneity and Universe Selection}
\label{subsec:endogeneity}

A central concern in dynamically constructed portfolios is whether changes in the investable
universe are influenced by future returns or by variables mechanically correlated with
subsequent performance. In the proposed framework, asset eligibility is determined
exclusively by observable variables that are fully lagged with respect to portfolio formation.

Eligibility depends on historical prices, trading volumes, and accounting fundamentals,
all observed at or before time $t-1$. No realized returns, forward returns, or ex post
performance measures enter the eligibility rule. As a consequence, the eligible universe
$\mathcal{U}_t$ is fully determined prior to portfolio optimization at time $t$.

To avoid implicit or undocumented endogeneity arising from interactions between universe
selection and factor exposure, the eligibility mechanism is explicitly formalized as a
deterministic function of an observable market state vector (Section~\ref{sec:formalization}).
This representation makes the conditioning structure transparent and separates state
determination from portfolio optimization.

\subsection{Look-Ahead Bias}
\label{subsec:lookahead}

All portfolio decisions are based exclusively on information available at the time of
rebalancing. Factor signals and eligibility criteria are computed using data observed at the
close of trading day $t-1$, while portfolio weights are implemented at time $t$.

Table~\ref{tab:timing} summarizes the information flow and timing assumptions adopted
throughout the framework.

\subsection{Data-Snooping and Parameter Choice}
\label{subsec:datasnooping}

The framework does not rely on in-sample optimization of portfolio weights or on
data-driven parameter tuning. Eligibility thresholds, factor definitions, and rebalancing
frequency are selected ex ante based on standard values commonly adopted in the
asset pricing and portfolio construction literature.

To further mitigate data-snooping concerns, empirical evaluation includes the deflated
Sharpe ratio, which adjusts for multiple testing and model selection effects. This statistic
provides a conservative assessment of risk-adjusted performance and guards against
overstating statistical significance in the presence of multiple design choices. 

\begin{table}[H]
	\centering
	\caption{Information timing and usage}
	\label{tab:timing}
	\begin{tabular}{lll}
		\toprule
		Quantity & Observed at & Used for \\
		\midrule
		Prices & $t-1$ close & Factor signal construction \\
		Trading volume (ADV) & $t-1$ & Eligibility and liquidity constraints \\
		Fundamentals & Last reported & Screening and factor inputs \\
		Portfolio weights & $t$ & Execution and portfolio rebalancing \\
		\bottomrule
	\end{tabular}
\end{table}

\subsection{Factor Redundancy}
\label{subsec:redundancy}

The factor set employed—momentum, value, and quality—consists of widely studied
characteristics with established economic interpretations. Signals are standardized
cross-sectionally at each rebalance date to ensure comparability and to prevent dominance
by any single factor.

To assess redundancy, we report the cross-factor correlation matrix and the marginal
contribution of each factor to the portfolio Sharpe ratio. These diagnostics allow us to
verify that portfolio behavior is not driven by a single latent exposure and that factor
contributions remain complementary rather than redundant.

\section{Dynamic Eligibility as a State-Dependent Constraint}
\label{sec:formalization}

This study proposes a portfolio construction framework in which factor exposure is
explicitly conditioned on observable market states through a dynamic eligibility mechanism.
Rather than treating the investable universe as exogenous, universe membership is modeled
as a state-dependent constraint that shapes the domain over which factor signals operate. Consistent with the information timing summarized in Table~\ref{tab:timing}, all state
variables and eligibility indicators entering the constraint $\mathcal{U}_t$ are observed
at or before time $t-1$, ensuring that the feasible asset set is fully determined prior to
portfolio optimization at time $t$. The resulting construction can be interpreted as a form of state-dependent constrained optimization, positioned between static multi-factor models and fully adaptive allocation schemes. The interpretation of eligibility as a state-dependent constraint aligns with
the literature on adaptive markets and state-contingent portfolio choice, in
which portfolio feasibility and exposure depend on observable market conditions
rather than static assumptions \cite{lo2004, ang2002, brandt2010}.

	Let $\mathcal{U}_t$ denote the set of assets eligible for investment at time $t$, and let $S_t$ represent a vector of observable market state variables, including liquidity, volatility, and market breadth. For each asset $i$, let $f_i(t)$ denote the vector of factor signals observed at time $t$.
	
	Eligibility is defined as
	\begin{equation}
		\mathcal{U}_t = \left\{ i \; : \; g_i(S_{t-1}) = 1 \right\},
	\end{equation}
	where $g_i(\cdot)$ is a deterministic indicator function based on lagged state variables.
	
	Portfolio weights are then obtained by solving
	\begin{equation}
		\mathbf{w}_t = \arg\max_{\mathbf{w}} \; \mathbf{w}^\top \mathbf{f}(t-1)
		\quad \text{s.t.} \quad
		\sum_i |w_i| = 1,\;\; i \in \mathcal{U}_t.
	\end{equation}
	
Under this formulation, factor exposure is conditioned on the prevailing market state
through the eligibility constraint $\mathcal{U}_t$. Changes in factor exposure arise
endogenously from changes in the feasible asset set rather than from time-varying factor
definitions or regime-specific parameter tuning.

This perspective places the proposed framework within the class of state-dependent
optimization and adaptive portfolio selection models, while preserving deterministic
rebalancing, bounded exposures, and full algorithmic transparency.

\subsection{Choice of State Variables}
\label{subsec:statevariables}

The selection of state variables entering the eligibility constraint is guided by
considerations of observability, operational relevance, and feasibility rather than
by optimality with respect to return prediction. The variables included in the state
vector—such as liquidity measures, volatility proxies, and cross-sectional breadth—
are chosen because they represent constraints that directly affect the implementability
and stability of real-world portfolios.

Importantly, these variables are not introduced as predictors of expected returns.
Instead, they proxy binding conditions under which portfolio construction must operate,
including trading capacity, turnover control, and exposure concentration. Their role is
therefore structural rather than predictive: they delimit the feasible asset set over
which factor signals are allowed to operate, without altering the factor definitions
themselves \cite{ferson1999, brandt2010}.

This design choice reflects a deliberate separation between portfolio feasibility and
return estimation. By conditioning eligibility on observable and lagged state variables,
the framework adapts endogenously to changing market conditions while avoiding reliance
on forecast-driven optimization or regime classification. Such constraint-based
constructions have been shown to improve stability and robustness in the presence of
estimation error and non-stationarity \cite{fabozzi2007, jagannathan2003}.

	\section{Counterfactual Baselines}
	\label{sec:baselines}
	
	To isolate the contribution of dynamic eligibility, we compare the proposed strategy against three counterfactual portfolios constructed using the same data and rebalancing schedule:
	
	\begin{enumerate}
		\item A fixed-universe multi-factor portfolio employing identical factor definitions but no eligibility constraints.
		\item An equal-weight portfolio, implemented both with and without dynamic eligibility.
		\item A market capitalization-weighted benchmark drawn from the same asset pool.
	\end{enumerate}
	
	All benchmarks use identical price data, rebalance dates, and transaction cost assumptions. Consequently, any performance differences can be attributed solely to the presence or absence of the dynamic eligibility mechanism.
	
	Results are summarized using one comparative performance figure and one table reporting risk-adjusted metrics.
	
	\section{Statistical Validation}
	\label{sec:statistics}
	
	Performance evaluation relies on a combination of standard and robustness-oriented statistics. Specifically, we report:
	
	\begin{itemize}
		\item Newey--West heteroskedasticity- and autocorrelation-consistent $t$-statistics for average returns;
		\item Deflated Sharpe ratios to account for multiple testing and model selection effects;
		\item Turnover-adjusted alpha, defined as excess return per unit of trading activity.
	\end{itemize}
	
	Together, these measures provide a conservative and methodologically robust assessment of economic significance. To control for multiple testing and selection effects, performance evaluation
	includes the deflated Sharpe ratio, which provides a conservative assessment
	of risk-adjusted returns under strategy selection uncertainty
	\cite{bailey2014}.

\section{Failure Modes and Regime Dependence}
\label{sec:failures}

The proposed framework does not dominate across all market environments.
Underperformance is observed in prolonged low-volatility bull markets, periods of
abundant liquidity, and phases characterized by pronounced factor crowding.

In such regimes, eligibility constraints become weakly binding and cross-sectional
dispersion diminishes, reducing the stabilizing effect of dynamic universe selection.
As a result, portfolio behavior converges toward that of conventional multi-factor or
equal-weight allocations, diminishing relative performance.

Documenting these failure modes clarifies the conditions under which dynamic eligibility
is most effective and reinforces the interpretation of the framework as a robustness-oriented
construction rather than an alpha-maximization strategy.

	\section{Investable Universe and Dynamic Inclusion}
	
	\subsection{Base universe}
	
	Let $\mathcal{U} = \{1,\dots,N\}$ denote a fixed reference universe of tradable assets, which may include equities, ADRs, ETFs, or other liquid instruments.
	
	\subsection{Dynamic inclusion rule}
	
	At each rebalance date $t$, the eligible universe $\mathcal{U}_t \subseteq \mathcal{U}$ is defined as:
	\[
	\mathcal{U}_t =
	\left\{
	i \in \mathcal{U} :
	\begin{aligned}
		& H_i(t) \ge H_{\min}, \\
		& \overline{ADV}_i(t; L_{ADV}) \ge ADV_{\min}
	\end{aligned}
	\right\},
	\]
	where $H_i(t)$ denotes available price history and $\overline{ADV}_i$ denotes average dollar volume over a fixed lookback window.
	
	Dynamic inclusion prevents systematic overweighting of immature or illiquid assets and mitigates survivorship and look-ahead bias\cite{brown1992,carhart1997,chen2002}. Liquidity constraints play a first-order role in determining feasible
	portfolio allocations and realized performance, particularly outside
	large-cap universes \cite{pastor2003, korajczyk2008}.

	\section{Baseline Weighting}
	
	For each rebalance date $t$, baseline weights are defined as equal weight over the eligible universe:
	\[
	w_i^{EW}(t) =
	\begin{cases}
		\frac{1}{|\mathcal{U}_t|}, & i \in \mathcal{U}_t, \\
		0, & i \notin \mathcal{U}_t.
	\end{cases}
	\]
	
	Equal weighting minimizes concentration risk and avoids dependence on market capitalization or estimated risk\cite{demiguel2009,plyakha2012,fernholz2002}.
	
	\section{Factor Signals}
	
	Let $\mathcal{F}$ denote a fixed set of factors, typically including momentum, value, and quality\cite{jegadeesh1993,fama1992,fama2015,novy2013}.
	
	For each factor $f \in \mathcal{F}$, a raw signal $s_{i,f}(t)$ is computed for each eligible asset. Signals are standardized cross-sectionally:
	\[
	z_{i,f}(t) = \frac{s_{i,f}(t) - \mu_f(t)}{\sigma_f(t)},
	\]
	where $\mu_f(t)$ and $\sigma_f(t)$ are cross-sectional moments over $\mathcal{U}_t$.
	
	\section{Bounded Multi-Factor Tilts}
	
	\subsection{Composite factor score}
	
	A convex combination of factor scores is formed:
	\[
	z_i(t) = \sum_{f \in \mathcal{F}} \alpha_f z_{i,f}(t),
	\quad
	\alpha_f \ge 0,
	\quad
	\sum_f \alpha_f = 1.
	\]
	
	\subsection{Bounded multiplicative adjustment}
	
	Factor information is applied multiplicatively to baseline weights:
	\[
	m_i(t) = \mathrm{clip}\left(1 + \lambda z_i(t), m_{\min}, m_{\max}\right),
	\]
	with fixed bounds $0 < m_{\min} \le 1 \le m_{\max}$.
	
	Final portfolio weights are:
	\[
	w_i(t) =
	\frac{w_i^{EW}(t)\, m_i(t)}
	{\sum_{j \in \mathcal{U}_t} w_j^{EW}(t)\, m_j(t)}.
	\]
	
	This construction ensures that factor signals influence allocations without allowing any asset to dominate the portfolio\cite{jagannathan2003,clarke2006,chow2011}. Imposing explicit bounds on portfolio weights and factor adjustments can
	improve out-of-sample stability by mitigating estimation error and limiting
	concentration risk, even when such constraints are not ex ante optimal
	\cite{boyd2004, jagannathan2003}.

	\section{Rebalancing Protocol}
	
	Rebalancing occurs on a fixed semi-annual schedule mapped deterministically to trading days. Between rebalance dates, weights are held constant.
	
	Semi-annual rebalancing balances responsiveness to structural changes against turnover, transaction costs, and tax considerations\cite{arnott2005,daryanani2008,bodie2018}.
	
	\section{Portfolio Evolution}
	
	Let $r_i(t)$ denote the single-period return of asset $i$. Portfolio returns are computed as:
	\[
	r_p(t) = \sum_{i \in \mathcal{U}} w_i(t^{-}) r_i(t),
	\]
	with optional linear transaction costs applied at rebalance dates as a function of portfolio turnover.
	
	\section{Extension to Smaller-Cap Universes}
	
	\subsection{Motivation}
	
	Smaller-cap universes present additional challenges, including thinner liquidity, higher volatility, noisier fundamentals, and increased transaction costs\cite{amihud2002,frazzini2015}. These features exacerbate the instability of mean--variance optimization and can materially erode gross factor premia once implementation frictions are accounted for\cite{asness2018}.
	
	\subsection{Parameter adaptation}
	
	The proposed framework extends naturally to smaller-cap universes by adjusting constraint parameters rather than altering structure. Typical adaptations include:
	
	\begin{itemize}
		\item higher minimum average dollar volume thresholds,
		\item tighter bounds on weight multipliers,
		\item reduced tilt strength $\lambda$,
		\item unchanged or lower rebalancing frequency.
	\end{itemize}
	
	No additional estimation layers are introduced.
	
	\subsection{Robustness considerations}
	
	Because the method relies on cross-sectional ranks, bounded adjustments, and deterministic rebalancing, it remains stable even when factor signals are noisy or transient, as is common in smaller-cap settings.
	
	\section{Algorithms}
	\subsection{Dynamic Inclusion and Bounded Multi-Factor Tilts}
	\label{subsec:dynamic_multifactor}
	
	This subsection introduces the core portfolio construction algorithm used throughout the paper. 
	The method combines dynamic asset inclusion, equal-weight baselines, and bounded multi-factor tilts under a fixed semi-annual rebalancing schedule. 
	The design explicitly avoids estimation of expected returns and covariance matrices, relying instead on cross-sectional factor rankings and hard weight constraints. The full procedure is defined by Algorithm 1a (dynamic inclusion and signal construction) and Algorithm 1b (weighting and rebalancing).
	
	\paragraph{Inputs.}
	The algorithm requires the following inputs:
	\begin{itemize}
		\item a reference asset universe $\mathcal{U}$;
		\item total-return price series $P_{t,i}$ for all $i \in \mathcal{U}$;
		\item (optional) trading volumes $V_{t,i}$ and fundamental data $F_{t,i}$;
		\item a deterministic set of rebalance dates $\mathcal{R}$;
		\item inclusion parameters: minimum history $H_{\min}$, minimum average dollar volume $ADV_{\min}$, and liquidity lookback $L_{ADV}$;
		\item factor parameters: factor set $\mathcal{F}$, mixture weights $\alpha_f$, tilt strength $\lambda$, and multiplier bounds $m_{\min}, m_{\max}$.
	\end{itemize}
	
	\paragraph{Outputs.}
	The algorithm produces:
	\begin{itemize}
		\item target portfolio weights $w_i(t)$ at each rebalance date $t \in \mathcal{R}$;
		\item a fully specified, piecewise-constant weight trajectory suitable for daily portfolio return computation.
	\end{itemize}
	
	Dynamic inclusion ensures that only assets satisfying liquidity and maturity criteria are eligible at each rebalance date. 
	Factor information influences allocations only through bounded multiplicative adjustments around an equal-weight baseline, preventing concentration and ensuring robustness under signal noise.
	The fixed semi-annual rebalancing schedule limits turnover and improves operational and tax efficiency.
	
\begin{algorithm}[H]
	\caption{Dynamic Inclusion and Factor Signal Construction}
	\label{alg:dynamic_inclusion}
	\begin{algorithmic}[1]
		
		\Require
		Asset universe $\mathcal{U}$, trading days $[T_0,T_1]$ \\
		Rebalance dates $\mathcal{R}$ \\
		Prices $P_{t,i}$, volumes $V_{t,i}$, fundamentals $F_{t,i}$ \\
		Inclusion parameters $H_{\min}, ADV_{\min}, L_{ADV}$
		
		\Ensure
		Eligible universe $\mathcal{U}_t$ and factor z-scores $z_{i,f}(t)$
		
		\For{each rebalance date $t \in \mathcal{R}$}
		
		\State \textbf{Dynamic inclusion}
		\State $\mathcal{U}_t \gets 
		\{ i \in \mathcal{U} :
		H_i(t) \ge H_{\min} \land
		\overline{ADV}_i(t;L_{ADV}) \ge ADV_{\min} \}$
		
		\For{each factor $f \in \mathcal{F}$}
		\For{each $i \in \mathcal{U}_t$}
		\State Compute raw signal $s_{i,f}(t)$
		\EndFor
		\State Standardize cross-sectionally:
		\State $z_{i,f}(t) \gets
		\dfrac{s_{i,f}(t) - \mu_f(t)}{\sigma_f(t)}$
		\EndFor
		
		\EndFor
		
		\State \Return $\{\mathcal{U}_t\}$ and $\{z_{i,f}(t)\}$
		
	\end{algorithmic}
\end{algorithm}

\begin{algorithm}[H]
	\caption{Semi-Annual Rebalancing with Bounded Multi-Factor Tilts}
	\label{alg:weighting_rebalance}
	\begin{algorithmic}[1]
		
		\Require
		Eligible universe $\mathcal{U}_t$, factor scores $z_{i,f}(t)$ \\
		Factor weights $\alpha_f$, tilt strength $\lambda$ \\
		Bounds $m_{\min}, m_{\max}$
		
		\Ensure
		Portfolio weights $w_{t,i}$ and returns $r_p(t)$
		
		\For{each rebalance date $t \in \mathcal{R}$}
		
		\State \textbf{Equal-weight baseline}
		\For{each $i \in \mathcal{U}_t$}
		\State $w_i^{EW} \gets 1 / |\mathcal{U}_t|$
		\EndFor
		
		\State \textbf{Composite factor score}
		\For{each $i \in \mathcal{U}_t$}
		\State $z_i(t) \gets \sum_{f \in \mathcal{F}} \alpha_f z_{i,f}(t)$
		\EndFor
		
		\State \textbf{Bounded tilt}
		\For{each $i \in \mathcal{U}_t$}
		\State $m_i(t) \gets
		\mathrm{clip}\!\left(1+\lambda z_i(t), m_{\min}, m_{\max}\right)$
		\EndFor
		
		\State \textbf{Final weights}
		\State $S \gets \sum_{j \in \mathcal{U}_t} w_j^{EW} m_j(t)$
		\For{each $i \in \mathcal{U}$}
		\If{$i \in \mathcal{U}_t$}
		\State $w_i(t) \gets \dfrac{w_i^{EW} m_i(t)}{S}$
		\Else
		\State $w_i(t) \gets 0$
		\EndIf
		\EndFor
		
		\EndFor
		
		\State \Return $\{w_{t,i}\}$
		
	\end{algorithmic}

\end{algorithm}

\subsection{Small-Cap Extension with Liquidity-Weighted Caps}
\label{subsec:smallcap_extension}

This subsection presents an extension of the baseline algorithm tailored to smaller-cap universes, where liquidity constraints and transaction costs are more pronounced.
The core structure of dynamic inclusion and bounded factor tilts is preserved, while an additional liquidity-dependent cap is imposed on individual asset weights.

\paragraph{Inputs.}
In addition to the inputs required by the baseline algorithm, the small-cap variant requires:
\begin{itemize}
	\item a liquidity proxy for each eligible asset, typically average dollar volume $\overline{ADV}_i(t)$;
	\item cap parameters: a global maximum weight $c_{\max}$, a scaling coefficient $\kappa$, and an elasticity parameter $\gamma$.
\end{itemize}

\paragraph{Outputs.}
The algorithm outputs:
\begin{itemize}
	\item liquidity-feasible portfolio weights $w_i(t)$ that satisfy both factor-based bounds and liquidity-weighted caps.
\end{itemize}

	\begin{algorithm}[H]
		\caption{Small-Cap Variant with Liquidity-Weighted Weight Caps}
		\label{alg:smallcap_liqcap}
		\begin{algorithmic}[1]
			\Require
			Eligible universe $\mathcal{U}_t$ at each rebalance date $t \in \mathcal{R}$ \\
			Baseline weights $w_i^{EW}(t)$, factor z-scores $z_{i,f}(t)$ \\
			Factor weights $\alpha_f$, tilt strength $\lambda$ \\
			Multiplier bounds $m_{\min}, m_{\max}$ \\
			Liquidity proxy $\overline{ADV}_i(t)$ over lookback $L_{ADV}$ \\
			Cap parameters: $c_{\max} \in (0,1)$, $\kappa>0$, $\gamma \in [0,1]$, tolerance $\varepsilon$
			\Ensure
			Liquidity-feasible target weights $w_i(t)$
			
			\For{each rebalance date $t \in \mathcal{R}$}
			
			\State \textbf{Composite factor score}
			\For{each $i \in \mathcal{U}_t$}
			\State $z_i(t) \gets \sum_{f \in \mathcal{F}} \alpha_f z_{i,f}(t)$
			\State $m_i(t) \gets \mathrm{clip}\!\left(1+\lambda z_i(t),\, m_{\min},\, m_{\max}\right)$
			\State $w_i^{raw} \gets w_i^{EW}(t)\, m_i(t)$
			\EndFor
			\State Normalize: $w_i^{raw} \gets w_i^{raw}/\sum_{j\in\mathcal{U}_t} w_j^{raw}$
			
			\State \textbf{Liquidity-weighted caps}
			\State $ADV^{med}(t) \gets \mathrm{median}_{j\in\mathcal{U}_t}\ \overline{ADV}_j(t)$
			\For{each $i \in \mathcal{U}_t$}
			\State $c_i(t) \gets \min\!\left(c_{\max},\ \kappa \cdot \left(\frac{\overline{ADV}_i(t)}{ADV^{med}(t)}\right)^{\gamma}\right)$
			\State $c_i(t) \gets \min\!\left(c_i(t),\, 1\right)$
			\EndFor
			
			\State \textbf{Cap-and-redistribute (iterative projection)}
			\State Initialize $w_i \gets w_i^{raw}$ for all $i \in \mathcal{U}_t$
			\Repeat
			\State $\Delta \gets 0$
			\State $\mathcal{B} \gets \{ i \in \mathcal{U}_t : w_i > c_i(t) \}$ \Comment{Breaching set}
			\If{$\mathcal{B} = \emptyset$}
			\State \textbf{break}
			\EndIf
			\State $E \gets \sum_{i \in \mathcal{B}} (w_i - c_i(t))$ \Comment{Excess mass}
			\For{each $i \in \mathcal{B}$}
			\State $w_i \gets c_i(t)$
			\EndFor
			\State $\mathcal{F}_t \gets \mathcal{U}_t \setminus \mathcal{B}$ \Comment{Free set}
			\State $S \gets \sum_{j \in \mathcal{F}_t} w_j$
			\For{each $j \in \mathcal{F}_t$}
			\State $w_j \gets w_j + E \cdot \frac{w_j}{S}$ \Comment{Redistribute proportionally}
			\EndFor
			\State $\Delta \gets E$
			\Until{$\Delta < \varepsilon$}
			
			\State Set $w_i(t)\gets w_i$ for $i\in\mathcal{U}_t$ and $w_i(t)\gets 0$ for $i\notin\mathcal{U}_t$
			\EndFor
			
			\State \Return $\{w_i(t)\}$
		\end{algorithmic}
	\end{algorithm}
	
	\subsection{Optional IC/IR Diagnostics for Factor Signals}
	\label{subsec:icir_calibration}
	
	This subsection describes an optional diagnostic procedure for assessing the relative
	cross-sectional predictive strength of factor signals. The procedure is not required for
	portfolio construction and does not alter the core algorithm presented in Sections 8--12.
	
	The method is designed to operate independently of the portfolio construction algorithm and may be applied periodically or used as a robustness check. Factor mixture weights are calibrated in two stages: construction of information coefficient time series from forward returns (Algorithm A3a), followed by information-ratio aggregation and mapping to non-negative factor weights (Algorithm A3b).
	
	\paragraph{Inputs.}
	The calibration algorithm requires:
	\begin{itemize}
		\item rebalance dates $\mathcal{R}$ and eligible universes $\mathcal{U}_t$;
		\item cross-sectional factor z-scores $z_{i,f}(t)$ for each factor $f \in \mathcal{F}$;
		\item total-return prices $P_{t,i}$ to compute forward returns over a fixed horizon $H$;
		\item a minimum observation threshold $M_{\min}$ to ensure statistical stability.
	\end{itemize}
	
	\paragraph{Outputs.}
	The algorithm returns:
	\begin{itemize}
		\item non-negative factor mixture weights $\alpha_f$ satisfying $\sum_f \alpha_f = 1$.
	\end{itemize}
	
	\begin{algorithm}[H]
		\caption{Forward Return Computation and Information Coefficient Construction}
		\label{alg:ic_construction}
		\begin{algorithmic}[1]
			
			\Require
			Rebalance dates $\mathcal{R}$, eligible sets $\mathcal{U}_t$ \\
			Factor z-scores $z_{i,f}(t)$ for $f \in \mathcal{F}$, $i \in \mathcal{U}_t$ \\
			Forward return horizon $H$ (trading days) \\
			Total return prices $P_{t,i}$
			
			\Ensure
			Time series of information coefficients $\mathcal{I}_f$ for each factor $f \in \mathcal{F}$
			
			\vspace{0.5em}
			\State \textbf{Compute forward returns at rebalance dates}
			\For{each $t \in \mathcal{R}$}
			\For{each $i \in \mathcal{U}_t$}
			\State $R_i^{fwd}(t) \gets \frac{P_{t+H,i}}{P_{t,i}} - 1$
			\EndFor
			\EndFor
			
			\vspace{0.5em}
			\State \textbf{Compute information coefficients}
			\For{each factor $f \in \mathcal{F}$}
			\State Initialize list $\mathcal{I}_f \gets [\,]$
			\For{each $t \in \mathcal{R}$}
			\State Let $x_i \gets z_{i,f}(t)$ and $y_i \gets R_i^{fwd}(t)$ for $i \in \mathcal{U}_t$
			\If{$|\mathcal{U}_t|$ is sufficient and $(x_i, y_i)$ not missing}
			\State $IC_f(t) \gets \mathrm{SpearmanCorr}\big(\{x_i\}, \{y_i\}\big)$
			\State Append $IC_f(t)$ to $\mathcal{I}_f$
			\EndIf
			\EndFor
			\EndFor
			
			\State \Return $\{\mathcal{I}_f\}_{f \in \mathcal{F}}$
			
		\end{algorithmic}
	\end{algorithm}
	
	\begin{algorithm}[H]
		\caption{Information Ratio Aggregation and Mapping to Factor Weights}
		\label{alg:ir_weight_mapping}
		\begin{algorithmic}[1]
			
			\Require
			Information coefficient time series $\mathcal{I}_f$ for $f \in \mathcal{F}$ \\
			Minimum observation threshold $M_{\min}$
			
			\Ensure
			Factor mixture weights $\alpha_f$ with $\alpha_f \ge 0$ and $\sum_f \alpha_f = 1$
			
			\vspace{0.5em}
			\State \textbf{Compute information ratios}
			\For{each factor $f \in \mathcal{F}$}
			\If{$|\mathcal{I}_f| < M_{\min}$ \textbf{or} $\mathrm{Std}(\mathcal{I}_f)=0$}
			\State $IR_f \gets 0$
			\Else
			\State $IR_f \gets
			\dfrac{\mathrm{Mean}(\mathcal{I}_f)}{\mathrm{Std}(\mathcal{I}_f)}$
			\EndIf
			\EndFor
			
			\vspace{0.5em}
			\State \textbf{Map IR to non-negative mixture weights}
			\For{each factor $f \in \mathcal{F}$}
			\State $s_f \gets \max(IR_f, 0)$
			\EndFor
			
			\If{$\sum_{f \in \mathcal{F}} s_f = 0$}
			\State $\alpha_f \gets 1 / |\mathcal{F}|$ for all $f \in \mathcal{F}$
			\Else
			\For{each factor $f \in \mathcal{F}$}
			\State $\alpha_f \gets \dfrac{s_f}{\sum_{g \in \mathcal{F}} s_g}$
			\EndFor
			\EndIf
			
			\State \Return $\{\alpha_f\}_{f \in \mathcal{F}}$
			
		\end{algorithmic}
	\end{algorithm}
	
All empirical results reported in this paper are obtained without reliance on this calibration
procedure, which is included solely as a robustness and interpretability aid.

\subsection{Construction of Momentum, Value, and Quality Factors}
\label{subsec:factor_construction}

This subsection documents the construction of the factor signals used in the portfolio weighting algorithm (Algorithm~\ref{alg:factor_construction}). 
The objective is to obtain robust, cross-sectionally comparable signals that can be applied uniformly across assets, markets, and capitalization segments without introducing additional estimation layers. Factor construction proceeds in two stages: raw signal computation (Algorithm A2a) followed by cross-sectional cleaning and standardization (Algorithm A2b).

\paragraph{Data requirements.}
The factor construction requires the following inputs:
\begin{itemize}
	\item daily total-return price series for all assets in the eligible universe;
	\item market capitalization estimates at rebalance dates;
	\item fundamental accounting data reported at quarterly or annual frequency, including balance-sheet and income-statement items;
	\item a deterministic mapping between rebalance dates and the most recent available fundamental reports.
\end{itemize}

All data are used on an as-of basis to avoid look-ahead bias.

\paragraph{Momentum.}
Momentum is computed as trailing total return over a fixed lookback window with a short skip period,
\[
s_{i,\mathrm{MOM}}(t) = \frac{P_{t-S,i}}{P_{t-L_{\mathrm{mom}}-S,i}} - 1,
\]
corresponding to the standard 12--1 momentum definition used in the empirical literature.
The skip period mitigates short-term mean reversion and microstructure effects.
Momentum is computed exclusively from price data and is therefore available for all assets with sufficient history.

\paragraph{Value.}
Value is measured using the book-to-market ratio,
\[
s_{i,\mathrm{VAL}}(t) = \frac{\mathrm{BookEquity}_i(q)}{\mathrm{MktCap}_i(t)},
\]
where $q \le t$ denotes the most recent reporting date not older than a fixed staleness threshold.
Book-to-market is preferred over earnings-based measures in a global setting because it is less sensitive to negative earnings and cross-country accounting heterogeneity.
If no valid fundamentals are available within the staleness window, the value signal is treated as missing.

\paragraph{Quality.}
Quality is constructed as a composite of profitability and balance-sheet strength metrics,
\[
s_{i,\mathrm{QUAL}}(t) =
z(\mathrm{ROE}_i(q)) +
z(\mathrm{GrossMargin}_i(q)) +
z(-\mathrm{DebtToAssets}_i(q)).
\]
This formulation captures both operating efficiency and financial conservatism.
Each component is standardized cross-sectionally before aggregation to prevent any single metric from dominating the composite.

\paragraph{Cross-sectional normalization and robustness.}
Raw factor signals are winsorized at fixed quantile thresholds to reduce the influence of outliers and accounting anomalies.
For each factor, standardized z-scores are computed across the eligible universe at each rebalance date.
If a factor exhibits zero cross-sectional variance, its z-scores are set to zero, corresponding to a neutral exposure.

\paragraph{Missing data and neutral treatment.}
Assets with missing or stale fundamentals are not excluded from the portfolio.
Instead, missing factor signals result in neutral z-scores, ensuring that such assets retain their baseline equal weight and preventing unintended concentration effects.

\paragraph{Applicability to small-cap universes.}
In smaller-cap universes, factor signals may exhibit higher noise and lower persistence.
The proposed construction remains applicable, provided that stricter liquidity filters, stronger winsorization, and tighter weight bounds are applied at the portfolio-construction stage.
No changes to the factor definitions themselves are required.

\paragraph{Implementation considerations.}
The factor construction is fully cross-sectional and non-parametric.
It does not rely on time-series estimation, covariance modeling, or predictive regressions.
As a result, the signals integrate naturally with bounded multi-factor tilts and dynamic inclusion rules without amplifying estimation error.

		\begin{algorithm}[H]
		\caption{Cross-Sectional Cleaning and Standardization of Factor Signals}
		\label{alg:factor_standardization}
		\begin{algorithmic}[1]
			
			\Require
			Raw factor signals $s_{i,f}(t)$ for $i \in \mathcal{U}_t$ \\
			Winsorization level $p$
			
			\Ensure
			Cross-sectional factor z-scores $z_{i,f}(t)$
			
			\vspace{0.5em}
			\For{each factor $f \in \{\mathrm{MOM}, \mathrm{VAL}, \mathrm{QUAL}\}$}
			
			\State Winsorize $\{s_{i,f}(t)\}_{i\in\mathcal{U}_t}$
			to quantiles $[p,1-p]$
			
			\State Compute cross-sectional mean $\mu_f(t)$
			and standard deviation $\sigma_f(t)$
			
			\For{each $i \in \mathcal{U}_t$}
			\If{$\sigma_f(t) > 0$}
			\State $z_{i,f}(t) \gets
			\dfrac{s_{i,f}(t) - \mu_f(t)}{\sigma_f(t)}$
			\Else
			\State $z_{i,f}(t) \gets 0$
			\EndIf
			\EndFor
			
			\EndFor
			
			\State \Return $\{z_{i,f}(t)\}_{i\in\mathcal{U}_t,\,f}$
			
		\end{algorithmic}
	\end{algorithm}
	
	\begin{algorithm}[H]
		\caption{Construction of Momentum, Value, and Quality Factor Signals}
		\label{alg:factor_construction}
		\begin{algorithmic}[1]
			
			\Require
			Eligible universe $\mathcal{U}_t$ at rebalance date $t$ \\
			Total-return prices $P_{d,i}$ for trading days $d$ \\
			Market capitalization $\mathrm{MktCap}_i(t)$ \\
			Fundamental data $F_{q,i}$ reported at dates $q$ \\
			Momentum parameters: lookback $L_{\mathrm{mom}}$, skip $S$ \\
			Fundamental staleness limit $L_{\mathrm{fund}}$ \\
			Winsorization level $p$
			
			\Ensure
			Cross-sectional factor z-scores $z_{i,f}(t)$ for
			$f \in \{\mathrm{MOM}, \mathrm{VAL}, \mathrm{QUAL}\}$
			
			\vspace{0.5em}
			\State \textbf{Momentum factor}
			\For{each $i \in \mathcal{U}_t$}
			\State $s_{i,\mathrm{MOM}}(t) \gets
			\dfrac{P_{t-S,i}}{P_{t-L_{\mathrm{mom}}-S,i}} - 1$
			\EndFor
			
			\vspace{0.5em}
			\State \textbf{Value factor}
			\For{each $i \in \mathcal{U}_t$}
			\State Select most recent fundamentals $F_{q,i}$ with
			$q \le t$ and $(t-q) \le L_{\mathrm{fund}}$
			\If{fundamentals available}
			\State $s_{i,\mathrm{VAL}}(t) \gets
			\dfrac{\mathrm{BookEquity}_i(q)}{\mathrm{MktCap}_i(t)}$
			\Else
			\State $s_{i,\mathrm{VAL}}(t) \gets \text{NaN}$
			\EndIf
			\EndFor
			
			\vspace{0.5em}
			\State \textbf{Quality factor}
			\For{each $i \in \mathcal{U}_t$}
			\If{required fundamentals available}
			\State $s_{i,\mathrm{QUAL}}(t) \gets
			z(\mathrm{ROE}_i(q))
			+ z(\mathrm{GrossMargin}_i(q))
			+ z(-\mathrm{DebtToAssets}_i(q))$
			\Else
			\State $s_{i,\mathrm{QUAL}}(t) \gets \text{NaN}$
			\EndIf
			\EndFor
			
			\vspace{0.5em}
			\State \textbf{Cross-sectional cleaning and standardization}
			\For{each factor $f \in \{\mathrm{MOM}, \mathrm{VAL}, \mathrm{QUAL}\}$}
			\State Winsorize $\{s_{i,f}(t)\}_{i\in\mathcal{U}_t}$ to quantiles $[p,1-p]$
			\State Compute mean $\mu_f(t)$ and standard deviation $\sigma_f(t)$
			\For{each $i \in \mathcal{U}_t$}
			\If{$\sigma_f(t) > 0$}
			\State $z_{i,f}(t) \gets \dfrac{s_{i,f}(t) - \mu_f(t)}{\sigma_f(t)}$
			\Else
			\State $z_{i,f}(t) \gets 0$
			\EndIf
			\EndFor
			\EndFor
			
			\State \Return $\{z_{i,f}(t)\}_{i\in\mathcal{U}_t,\,f}$
			
		\end{algorithmic}
	\end{algorithm}

	\FloatBarrier
	\section{Diagrams}
	
\begin{figure}[H]
	\centering
	\begin{tikzpicture}[
		scale=0.9,
		transform shape,
		font=\small,
		node distance=7mm and 12mm,
		box/.style={draw, rounded corners=2pt, align=left, inner sep=5pt, text width=9.0cm},
		term/.style={draw, rounded corners=2pt, align=center, inner sep=5pt, text width=4.0cm},
		arrow/.style={-Latex, thick}
		]
		
		\node[term] (start) {Start (rebalance date $t\in\mathcal{R}$)};
		\node[box, below=of start] (incl) {\textbf{Dynamic inclusion}\\
			Compute eligible set:\\
			$\mathcal{U}_t \gets \{ i\in\mathcal{U}:\ H_i(t)\ge H_{\min}\ \land\ \overline{ADV}_i(t;L_{ADV})\ge ADV_{\min}\}$};
		\node[box, below=of incl] (raw) {\textbf{For each factor $f\in\mathcal{F}$}\\
			For each $i\in\mathcal{U}_t$: compute raw signal $s_{i,f}(t)$};
		\node[box, below=of raw] (std) {\textbf{Cross-sectional standardization}\\
			Compute $\mu_f(t),\ \sigma_f(t)$ over $\mathcal{U}_t$ and set\\
			$z_{i,f}(t)\gets \dfrac{s_{i,f}(t)-\mu_f(t)}{\sigma_f(t)}$};
		\node[term, below=of std] (out) {Output: $\mathcal{U}_t$ and $\{z_{i,f}(t)\}$};
		
		\draw[arrow] (start) -- (incl);
		\draw[arrow] (incl) -- (raw);
		\draw[arrow] (raw) -- (std);
		\draw[arrow] (std) -- (out);
		
	\end{tikzpicture}
	\caption{Flow diagram for Algorithm~\ref{alg:dynamic_inclusion}: dynamic inclusion and factor signal construction.}
	\label{fig:diag_dynamic_inclusion}
\end{figure}
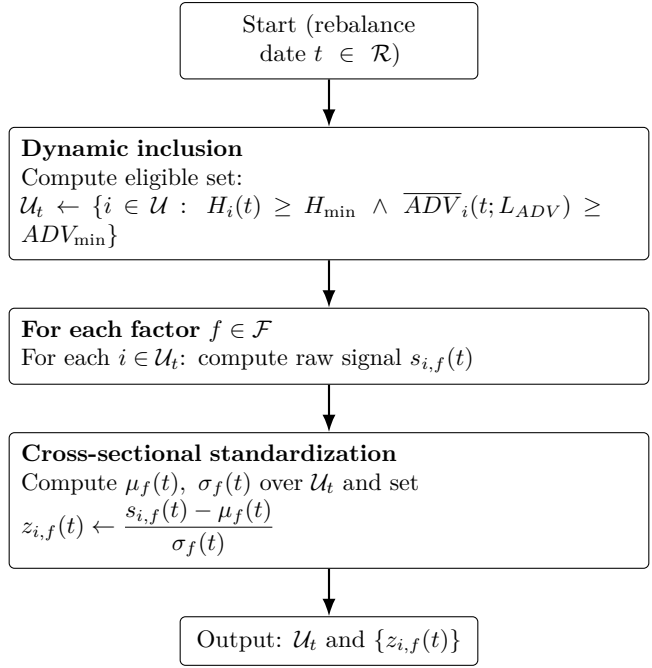
\begin{figure}[H]
	\centering
	\begin{tikzpicture}[
		scale=0.9,
		transform shape,
		font=\small,
		node distance=7mm and 12mm,
		box/.style={draw, rounded corners=2pt, align=left, inner sep=5pt, text width=9.0cm},
		term/.style={draw, rounded corners=2pt, align=center, inner sep=5pt, text width=4.5cm},
		arrow/.style={-Latex, thick}
		]
		
		\node[term] (start) {Start (rebalance date $t\in\mathcal{R}$)};
		\node[box, below=of start] (ew) {\textbf{Equal-weight baseline}\\
			For each $i\in\mathcal{U}_t$: $w_i^{EW}\gets 1/|\mathcal{U}_t|$};
		\node[box, below=of ew] (comp) {\textbf{Composite factor score}\\
			For each $i\in\mathcal{U}_t$:
			$z_i(t)\gets \sum_{f\in\mathcal{F}}\alpha_f\, z_{i,f}(t)$};
		\node[box, below=of comp] (tilt) {\textbf{Bounded multiplicative tilt}\\
			For each $i\in\mathcal{U}_t$:
			$m_i(t)\gets \mathrm{clip}(1+\lambda z_i(t),\, m_{\min},\, m_{\max})$};
		\node[box, below=of tilt] (norm) {\textbf{Normalize to final weights}\\
			$S\gets \sum_{j\in\mathcal{U}_t} w_j^{EW} m_j(t)$\\
			For $i\in\mathcal{U}_t$: $w_i(t)\gets \dfrac{w_i^{EW} m_i(t)}{S}$; else $w_i(t)\gets 0$};
		\node[term, below=of norm] (out) {Output: $\{w_i(t)\}$};
		
		\draw[arrow] (start) -- (ew);
		\draw[arrow] (ew) -- (comp);
		\draw[arrow] (comp) -- (tilt);
		\draw[arrow] (tilt) -- (norm);
		\draw[arrow] (norm) -- (out);
		
	\end{tikzpicture}
	\caption{Flow diagram for Algorithm~\ref{alg:weighting_rebalance}: semi-annual rebalancing with bounded multi-factor tilts.}
	\label{fig:diag_weighting_rebalance}
\end{figure}
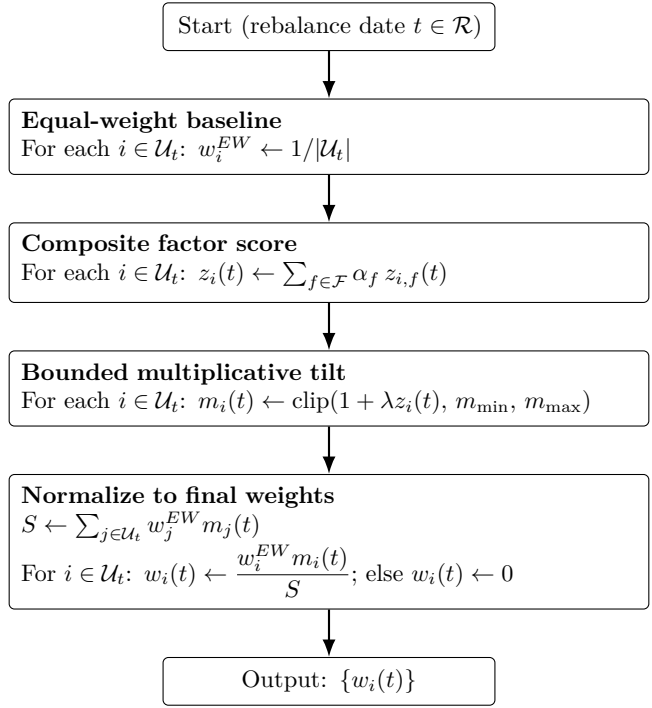

	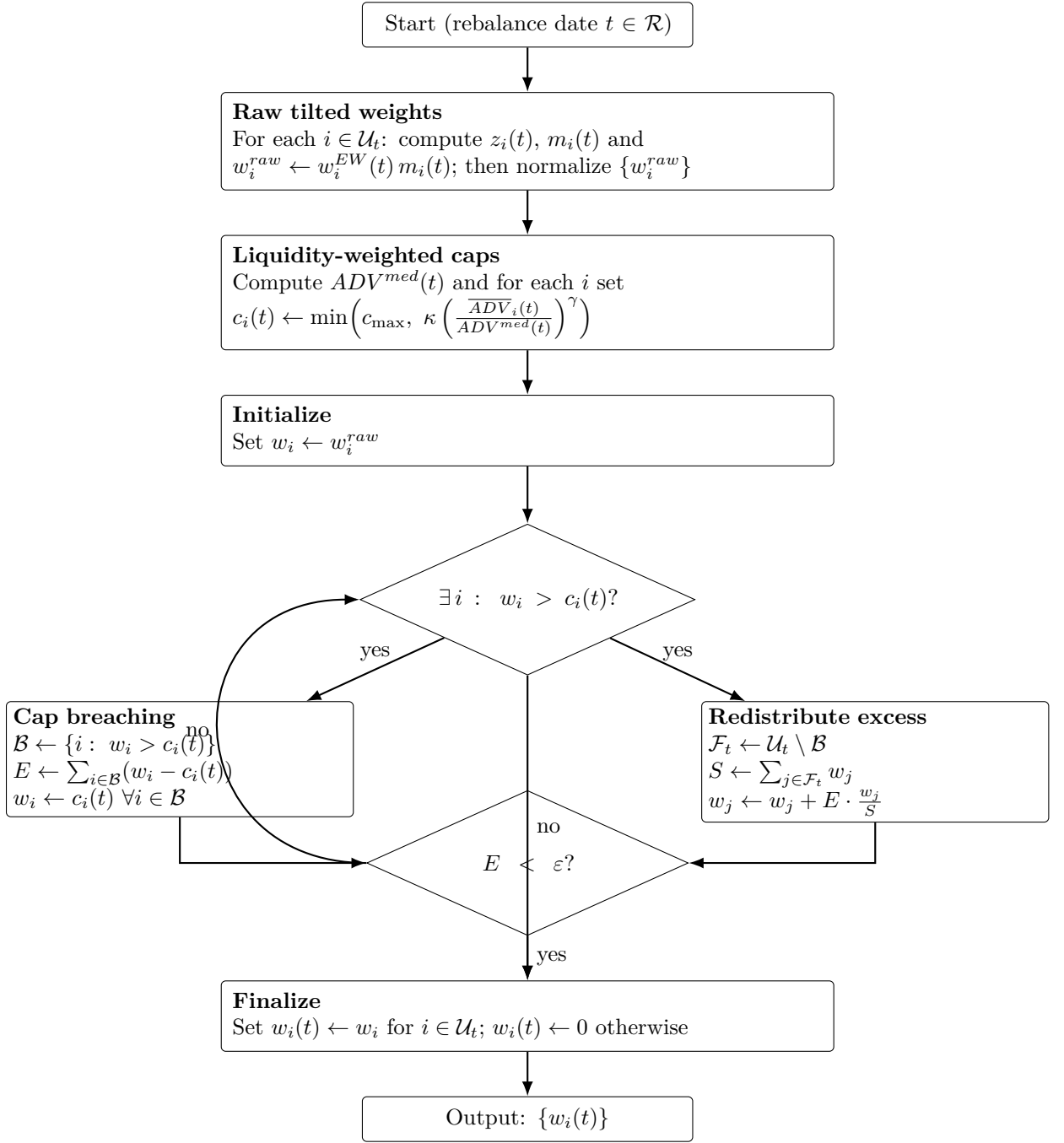
\begin{figure}[H]
		\centering
		\begin{tikzpicture}[
			font=\small,
			node distance=7mm and 12mm,
			box/.style={draw, rounded corners=2pt, align=left, inner sep=5pt, text width=9.2cm},
			smallbox/.style={draw, rounded corners=2pt, align=left, inner sep=3pt, text width=5.2cm},
			decision/.style={draw, diamond, aspect=2.2, align=center, inner sep=2pt, text width=4.0cm},
			term/.style={draw, rounded corners=2pt, align=center, inner sep=5pt, text width=4.8cm},
			arrow/.style={-Latex, thick}
			]
			\node[term] (start) {Start (rebalance date $t\in\mathcal{R}$)};
			\node[box, below=of start] (raw) {\textbf{Raw tilted weights}\\
				For each $i\in\mathcal{U}_t$: compute $z_i(t)$, $m_i(t)$ and\\
				$w_i^{raw}\gets w_i^{EW}(t)\, m_i(t)$; then normalize $\{w_i^{raw}\}$};
			\node[box, below=of raw] (caps) {\textbf{Liquidity-weighted caps}\\
				Compute $ADV^{med}(t)$ and for each $i$ set\\
				$c_i(t)\gets \min\!\left(c_{\max},\ \kappa\left(\frac{\overline{ADV}_i(t)}{ADV^{med}(t)}\right)^\gamma\right)$};
			\node[box, below=of caps] (init) {\textbf{Initialize}\\
				Set $w_i\gets w_i^{raw}$};
			\node[decision, below=9mm of init] (breach) {$\exists\, i:\ w_i>c_i(t)$?};
			
			\node[smallbox, below left=10mm and 14mm of breach] (cap) {%
				\textbf{Cap breaching}\\
				$\mathcal{B}\gets\{i:\ w_i>c_i(t)\}$\\
				$E\gets\sum_{i\in\mathcal{B}}(w_i-c_i(t))$\\
				$w_i\gets c_i(t)\ \forall i\in\mathcal{B}$
			};
			
			\node[smallbox, below right=10mm and 14mm of breach] (redis) {%
				\textbf{Redistribute excess}\\
				$\mathcal{F}_t\gets\mathcal{U}_t\setminus\mathcal{B}$\\
				$S\gets\sum_{j\in\mathcal{F}_t} w_j$\\
				$w_j\gets w_j + E\cdot \frac{w_j}{S}$
			};
			
			\node[decision, below=18mm of breach] (tol) {$E<\varepsilon$?};
			\node[box, below=of tol] (final) {\textbf{Finalize}\\
				Set $w_i(t)\gets w_i$ for $i\in\mathcal{U}_t$; $w_i(t)\gets 0$ otherwise};
			\node[term, below=of final] (out) {Output: $\{w_i(t)\}$};
			
			\draw[arrow] (start) -- (raw);
			\draw[arrow] (raw) -- (caps);
			\draw[arrow] (caps) -- (init);
			\draw[arrow] (init) -- (breach);
			
			\draw[arrow] (breach) -- node[midway, above] {yes} (cap);
			\draw[arrow] (breach) -- node[midway, above] {yes} (redis);
			
			\draw[arrow] (cap) |- (tol);
			\draw[arrow] (redis) |- (tol);
			
			\draw[arrow] (tol) -- node[midway, right] {yes} (final);
			\draw[arrow] (final) -- (out);
			
			% loop back when not within tolerance
			\draw[arrow] (tol.west) .. controls +(-3.0,0.0) and +(-3.0,0.0) .. node[midway, left] {no} (breach.west);
			
			% no-breach path
			\draw[arrow] (breach.south) -- node[midway, right] {no} (final.north);
		\end{tikzpicture}
		\caption{Flow diagram for Algorithm~\ref{alg:smallcap_liqcap}: small-cap variant with liquidity-weighted caps and iterative cap-and-redistribute projection.}
		\label{fig:diag_smallcap}
	\end{figure}

	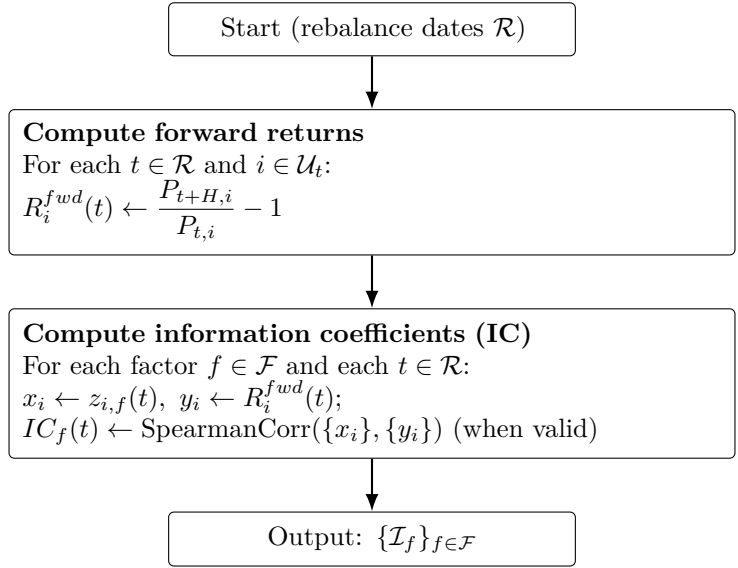
\begin{figure}[H]
		\centering
		\begin{tikzpicture}[
			font=\small,
			node distance=7mm and 12mm,
			box/.style={draw, rounded corners=2pt, align=left, inner sep=5pt, text width=9.2cm},
			term/.style={draw, rounded corners=2pt, align=center, inner sep=5pt, text width=5.0cm},
			arrow/.style={-Latex, thick}
			]
			\node[term] (start) {Start (rebalance dates $\mathcal{R}$)};
			\node[box, below=of start] (fwd) {\textbf{Compute forward returns}\\
				For each $t\in\mathcal{R}$ and $i\in\mathcal{U}_t$:\\
				$R_i^{fwd}(t)\gets \dfrac{P_{t+H,i}}{P_{t,i}}-1$};
			\node[box, below=of fwd] (ic) {\textbf{Compute information coefficients (IC)}\\
				For each factor $f\in\mathcal{F}$ and each $t\in\mathcal{R}$:\\
				$x_i\gets z_{i,f}(t),\ y_i\gets R_i^{fwd}(t)$;\\
				$IC_f(t)\gets \mathrm{SpearmanCorr}(\{x_i\},\{y_i\})$ (when valid)};
			\node[term, below=of ic] (out) {Output: $\{\mathcal{I}_f\}_{f\in\mathcal{F}}$};
			
			\draw[arrow] (start) -- (fwd);
			\draw[arrow] (fwd) -- (ic);
			\draw[arrow] (ic) -- (out);
		\end{tikzpicture}
		\caption{Flow diagram for Algorithm~\ref{alg:ic_construction}: forward-return computation and information coefficient construction.}
		\label{fig:diag_ic_construction}
	\end{figure}
	
	\begin{figure}[H]
		\centering
		\begin{tikzpicture}[
			font=\small,
			node distance=7mm and 12mm,
			box/.style={draw, rounded corners=2pt, align=left, inner sep=5pt, text width=9.2cm},
			smallbox/.style={draw, rounded corners=2pt, align=center, inner sep=3pt, text width=3.4cm},
			decision/.style={draw, diamond, aspect=3.2, align=center, inner sep=1pt, text width=2.4cm},
			term/.style={draw, rounded corners=2pt, align=center, inner sep=5pt, text width=5.2cm},
			arrow/.style={-Latex, thick}
			]
			
			\node[term] (start) {Start (IC series $\mathcal{I}_f$)};
			\node[box, below=of start] (ir) {\textbf{Compute information ratios}\\
				For each factor $f$: if $|\mathcal{I}_f|<M_{\min}$ or $\mathrm{Std}(\mathcal{I}_f)=0$, set
				$IR_f\gets 0$; else $IR_f\gets \dfrac{\mathrm{Mean}(\mathcal{I}_f)}{\mathrm{Std}(\mathcal{I}_f)}$};
			
			\node[box, below=of ir] (map) {\textbf{Map to non-negative scores}\\
				For each $f$: $s_f\gets \max(IR_f,0)$};
			
			\node[decision, below=10mm of map] (sum0) {$\sum_{f\in\mathcal{F}} s_f = 0$?};
			
			\node[smallbox, below left=7mm and -8mm of sum0] (eq) {%
				\textbf{Fallback}\\
				$\alpha_f\gets 1/|\mathcal{F}|$
			};
			
			\node[smallbox, below right=7mm and -8mm of sum0] (norm) {%
				\textbf{Normalize}\\
				$\alpha_f\gets \dfrac{s_f}{\sum_{g\in\mathcal{F}} s_g}$
			};
			
			\node[term, below=18mm of sum0] (out) {Output: $\{\alpha_f\}$};
			
			\draw[arrow] (start) -- (ir);
			\draw[arrow] (ir) -- (map);
			\draw[arrow] (map) -- (sum0);
			\draw[arrow] (sum0) -- node[midway, above] {yes} (eq);
			\draw[arrow] (sum0) -- node[midway, above] {no} (norm);
			\draw[arrow] (eq) |- (out);
			\draw[arrow] (norm) |- (out);
			
		\end{tikzpicture}
		\caption{Flow diagram for Algorithm~\ref{alg:ir_weight_mapping}: information-ratio aggregation and mapping to convex factor mixture weights.}
		\label{fig:diag_ir_mapping}
	\end{figure}
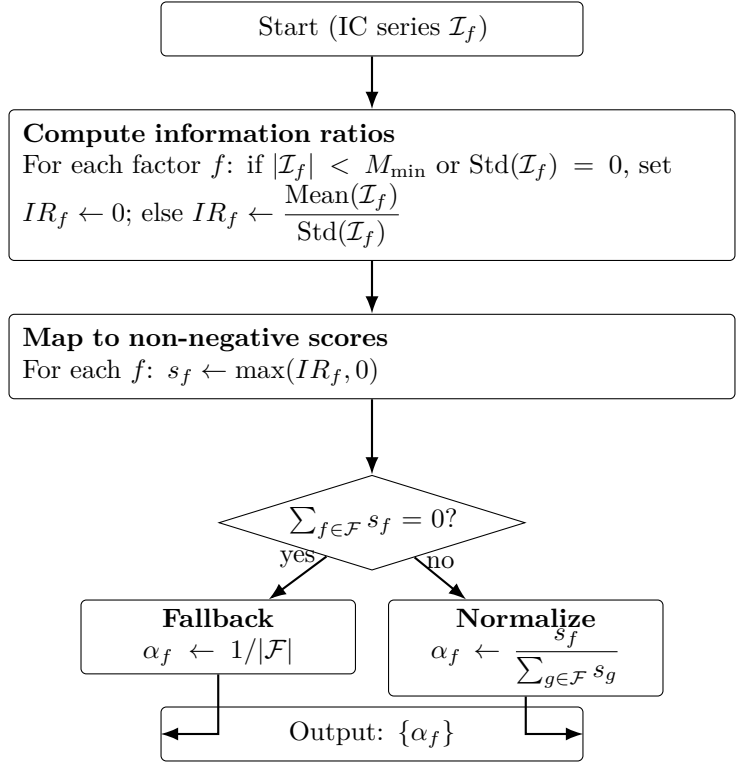
	
	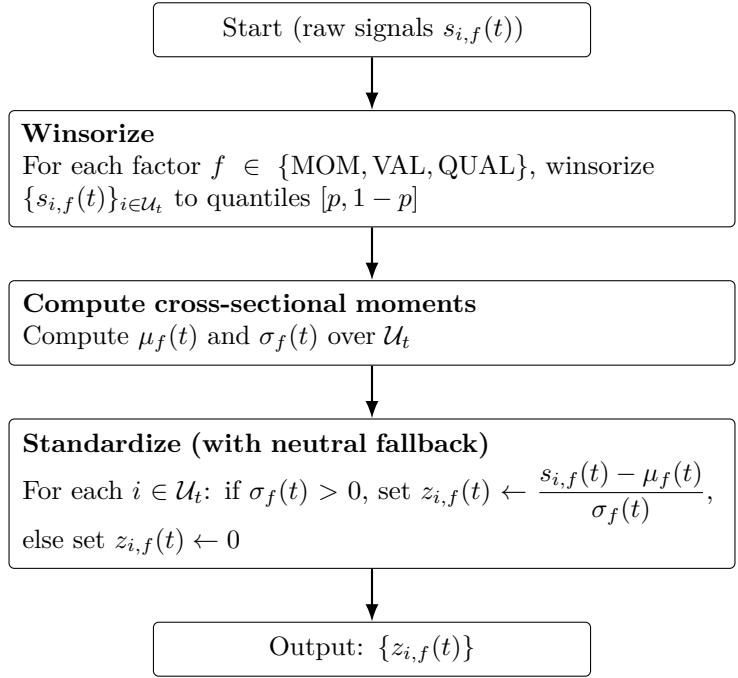
\begin{figure}[H]
		\centering
		\begin{tikzpicture}[
			font=\small,
			node distance=7mm and 12mm,
			box/.style={draw, rounded corners=2pt, align=left, inner sep=5pt, text width=9.2cm},
			term/.style={draw, rounded corners=2pt, align=center, inner sep=5pt, text width=5.4cm},
			arrow/.style={-Latex, thick}
			]
			\node[term] (start) {Start (raw signals $s_{i,f}(t)$)};
			\node[box, below=of start] (win) {\textbf{Winsorize}\\
				For each factor $f\in\{\mathrm{MOM},\mathrm{VAL},\mathrm{QUAL}\}$, winsorize
				$\{s_{i,f}(t)\}_{i\in\mathcal{U}_t}$ to quantiles $[p,1-p]$};
			\node[box, below=of win] (moments) {\textbf{Compute cross-sectional moments}\\
				Compute $\mu_f(t)$ and $\sigma_f(t)$ over $\mathcal{U}_t$};
			\node[box, below=of moments] (z) {\textbf{Standardize (with neutral fallback)}\\
				For each $i\in\mathcal{U}_t$:
				if $\sigma_f(t)>0$, set $z_{i,f}(t)\gets \dfrac{s_{i,f}(t)-\mu_f(t)}{\sigma_f(t)}$,
				else set $z_{i,f}(t)\gets 0$};
			\node[term, below=of z] (out) {Output: $\{z_{i,f}(t)\}$};
			
			\draw[arrow] (start) -- (win);
			\draw[arrow] (win) -- (moments);
			\draw[arrow] (moments) -- (z);
			\draw[arrow] (z) -- (out);
		\end{tikzpicture}
		\caption{Flow diagram for Algorithm~\ref{alg:factor_standardization}: cross-sectional cleaning and standardization of factor signals.}
		\label{fig:diag_factor_standardization}
	\end{figure}
	
	\FloatBarrier

	\section{Limitations}
	
	The framework does not target theoretical mean--variance efficiency and does not incorporate leverage or short selling. Results depend on factor definitions, data quality, and universe selection.
	
	\section{Scope of Applicability}
	
	The proposed framework is designed for robustness under estimation error,
	non-stationarity, and implementation constraints. As such, its suitability is
	inherently domain-specific. This section delineates the conditions under which
	the methodology is expected to perform as intended, as well as the settings in
	which alternative portfolio construction approaches are structurally more
	appropriate.
	
	\subsection{Market Regimes}
	
	The framework is well suited to environments characterized by moderate and
	persistent cross-sectional dispersion, where relative rankings across assets
	contain stable information. In such regimes, bounded multi-factor tilts allow
	the portfolio to benefit from systematic signals while limiting concentration
	risk and drawdown convexity.
	
	Conversely, the method is structurally disadvantaged in regimes dominated by a
	single, highly persistent factor or narrow leadership set (e.g.\ winner-take-most
	markets). In these settings, cap-weighted or unconstrained factor portfolios may
	exhibit superior upside capture. The proposed framework explicitly sacrifices
	such upside in favor of bounded exposure and robustness to subsequent regime
	reversals.
	
	In periods of severe market stress or policy-driven synchronization, when
	cross-sectional dispersion collapses, factor signals lose explanatory power.
	Under these conditions, the framework converges toward its equal-weight
	baseline, effectively reverting to a diversification-first allocation. This
	behavior is intentional and reflects the absence of exploitable relative
	information.
	
	\subsection{Universe Size and Structure}
	
	The methodology is most effective for universes of moderate to large size,
	where cross-sectional statistics are well-defined and liquidity constraints can
	be enforced without excessive turnover. Empirically, this corresponds to
	universes containing several dozen to several hundred liquid assets.
	
	For very small universes, cross-sectional normalization becomes unstable and
	factor rankings are overly sensitive to idiosyncratic effects. In such cases,
	the framework degenerates toward equal weighting with limited interpretability.
	At the opposite extreme, extremely large or heterogeneous universes lacking
	hierarchical segmentation may require additional structural layers (e.g.\
	regional or sectoral sleeves) to preserve robustness and implementability.
	
	\subsection{Asset Class Coverage}
	
	The framework is primarily intended for equity-like instruments, including
	single stocks, equity ETFs, and other ownership-based assets with comparable
	return distributions. While it may be adapted to certain commodity proxies, it
	is not designed for asset classes dominated by term-structure dynamics,
	convex payoffs, or leverage constraints, such as fixed income, foreign exchange,
	or derivative-intensive strategies.
	
	\subsection{Institutional and Operational Constraints}
	
	The method is compatible with long-horizon mandates that tolerate moderate
	tracking error and emphasize diversification, transparency, and turnover
	control. It is not suitable for benchmark-hugging mandates or strategies
	requiring tight ex-ante tracking-error control relative to capitalization-weighted
	indices.
	
	Similarly, the framework is not intended for high-frequency or tactical
	allocation contexts. Its deterministic rebalancing schedule and absence of
	volatility targeting or regime detection reflect a deliberate design choice in
	favor of stability over responsiveness.
	
	\subsection{Interpretation of Results}
	
	The framework should not be interpreted as an alpha-maximizing engine or a
	market-timing system. Its objective is to provide a robust core allocation that
	behaves predictably under adverse conditions, limits estimation risk, and
	remains operationally feasible across market regimes. Performance should
	therefore be evaluated in terms of concentration control, drawdown behavior,
	and robustness rather than short-term relative returns.
	
	\section{Core--Satellite Institutional Architecture Mapping}
	
	The proposed framework is naturally suited to serve as a \emph{core allocation}
	within a broader institutional portfolio architecture. In this role, it provides
	a stable, diversified foundation upon which more aggressive, tactical, or
	idiosyncratic strategies may be layered as satellites. This section clarifies how
	the methodology integrates with such architectures and delineates its functional
	boundaries relative to complementary portfolio components.
	
	\subsection{Role of the Core Allocation}
	
	Within a core--satellite structure, the core portfolio is responsible for
	delivering broad market exposure while controlling concentration, turnover, and
	estimation risk. The proposed framework fulfills this role through its use of an
	equal-weight baseline, bounded multi-factor tilts, and deterministic rebalancing.
	These features ensure that the core allocation remains interpretable, robust to
	regime changes, and operationally feasible over long horizons.
	
	Unlike return-maximizing or benchmark-tracking cores, the objective of the
	proposed method is not to dominate in all market environments, but to behave
	predictably under adverse conditions and to limit structural fragility arising
	from estimation error or excessive concentration.

	\subsection{Satellite Strategies and Complementarity}
	
	Satellite portfolios are intended to express views or exploit opportunities that
	the core framework deliberately avoids. These may include, but are not limited to,
	tactical momentum strategies, thematic or sector-focused allocations, macro-driven
	overlays, or idiosyncratic security selection.
	
	Such strategies typically accept higher concentration, higher turnover, or
	greater sensitivity to regime changes in exchange for potential upside. The
	presence of a robust core allocation allows these satellites to operate without
	dominating overall portfolio risk, as adverse outcomes in satellite components
	are partially absorbed by the stability of the core.

	\subsection{Interaction and Risk Allocation}
	
	The interaction between the core and satellite components is governed by explicit
	risk budgeting rather than signal integration. The core allocation establishes a
	stable risk baseline, while satellites are sized according to their expected
	volatility, drawdown potential, and correlation with the core portfolio.
	
	Importantly, the proposed framework does not require modification or internal
	state detection to accommodate satellite strategies. Its bounded and non-adaptive
	structure ensures that the core remains insulated from short-term tactical
	decisions, thereby preserving its role as a long-horizon allocator.
	
		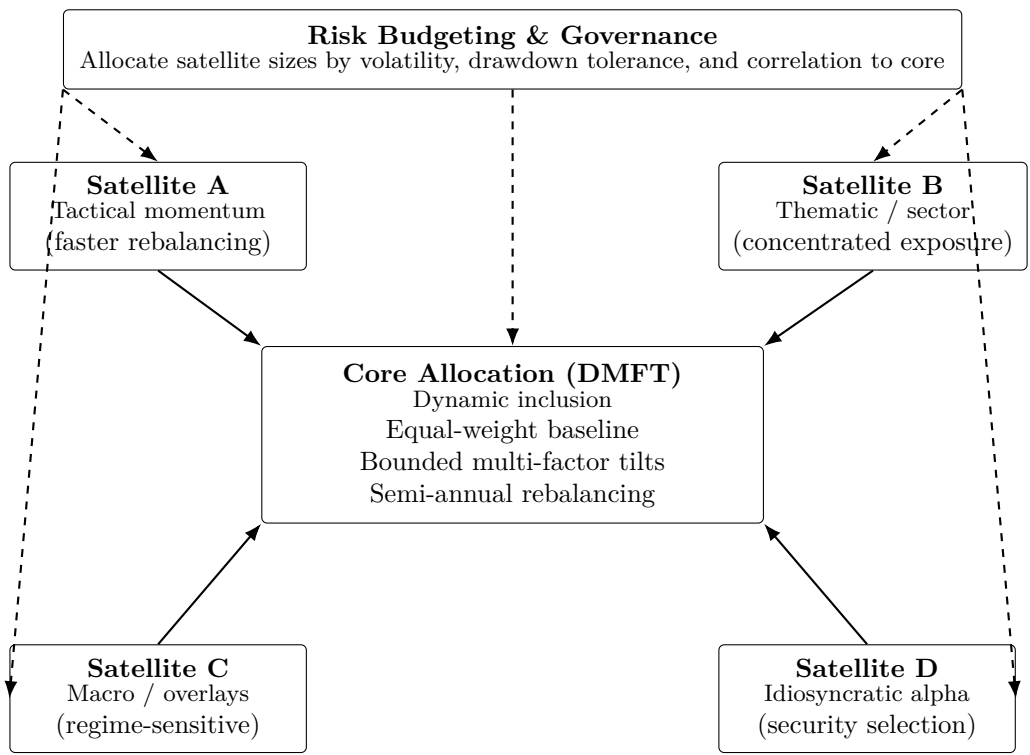
\begin{figure}[H]
		\centering
		\begin{tikzpicture}[
			font=\small,
			box/.style={draw, rounded corners=2pt, align=center, inner sep=6pt},
			sat/.style={draw, rounded corners=2pt, align=center, inner sep=5pt},
			arrow/.style={-Latex, thick},
			dasharrow/.style={-Latex, thick, dashed},
			node distance=14mm and 18mm
			]
			
			% --- Core (center) ---
			\node[box, minimum width=6.6cm, minimum height=2.25cm] (core) {%
				\textbf{Core Allocation (DMFT)}\\[-1mm]
				\footnotesize
				Dynamic inclusion\\
				Equal-weight baseline\\
				Bounded multi-factor tilts\\
				Semi-annual rebalancing
			};
			
			% --- Risk budgeting (top) ---
			% Put this sufficiently ABOVE satellites to avoid overlap with their boxes and labels.
			\node[box, above=34mm of core, minimum width=8.6cm] (budget) {%
				\textbf{Risk Budgeting \& Governance}\\[-1mm]
				\footnotesize Allocate satellite sizes by volatility, drawdown tolerance, and correlation to core
			};
			
			% --- Satellites (corners) ---
			% Place satellites below the budgeting box and away from the top margin.
			\node[sat, above left=10mm and -6mm of core, minimum width=3.9cm] (mom) {%
				\textbf{Satellite A}\\[-1mm]
				\footnotesize Tactical momentum\\
				(faster rebalancing)
			};
			
			\node[sat, above right=10mm and -6mm of core, minimum width=3.9cm] (them) {%
				\textbf{Satellite B}\\[-1mm]
				\footnotesize Thematic / sector\\
				(concentrated exposure)
			};
			
			% Give more vertical clearance at the bottom so nothing collides with the caption.
			\node[sat, below left=16mm and -6mm of core, minimum width=3.9cm] (macro) {%
				\textbf{Satellite C}\\[-1mm]
				\footnotesize Macro / overlays\\
				(regime-sensitive)
			};
			
			\node[sat, below right=16mm and -6mm of core, minimum width=3.9cm] (alpha) {%
				\textbf{Satellite D}\\[-1mm]
				\footnotesize Idiosyncratic alpha\\
				(security selection)
			};
			
			% --- Arrows: governance/risk-budget constraints (dashed) ---
			\draw[dasharrow] (budget.south west) -- (mom.north);
			\draw[dasharrow] (budget.south east) -- (them.north);
			\draw[dasharrow] (budget.south west) -- (macro.west);
			\draw[dasharrow] (budget.south east) -- (alpha.east);
			\draw[dasharrow] (budget.south) -- (core.north);
			
			% --- Arrows: satellites aggregate into core (solid) ---
			\draw[arrow] (mom.south) -- (core.north west);
			\draw[arrow] (them.south) -- (core.north east);
			\draw[arrow] (macro.north) -- (core.south west);
			\draw[arrow] (alpha.north) -- (core.south east);
			
			% NOTE: The in-figure annotation caused overlap in your output.
			% Keep the sentence in the caption instead (below) to eliminate collisions.
			
		\end{tikzpicture}
		
		\caption{Core--satellite institutional mapping of the proposed framework. The DMFT portfolio
			acts as a robust core allocation (dynamic inclusion, equal-weight baseline, bounded tilts,
			deterministic rebalancing). More aggressive or tactical strategies may be layered as satellites.
			Portfolio construction is governed by explicit risk budgeting (e.g., sizing by volatility,
			drawdown tolerance, and correlation to the core) rather than by integrating satellite signals
			into the core optimization. Solid arrows indicate conceptual aggregation toward the overall
			portfolio; dashed arrows indicate governance and risk-budget constraints. The core provides a
			stable risk baseline; satellites express views with higher concentration, turnover, or regime
			sensitivity.}
		\label{fig:core_satellite}
	\end{figure}
	
	\subsection{Interpretation and Institutional Relevance}
	
	Interpreted within a core--satellite architecture, the proposed framework should
	not be viewed as a standalone alpha engine, but as a structural component designed
	to anchor portfolio risk and diversification. Its value lies in providing a
	robust baseline that reduces the dependency of overall portfolio outcomes on the
	success of any single strategy or regime.
	
	This positioning is particularly relevant for institutional investors, wealth
	managers, and multi-strategy portfolios seeking a transparent and resilient core
	allocation that coexists coherently with faster or more aggressive investment
	approaches.

\section{Conclusion}

This paper has presented a fully specified portfolio construction framework designed to
prioritize robustness, interpretability, and operational feasibility under estimation error,
non-stationarity, and practical trading constraints. The methodology combines dynamic asset
inclusion, deterministic rebalancing, and bounded multi-factor tilts applied to an equal-weight
baseline, thereby avoiding explicit estimation of expected returns and covariance matrices.

A central contribution of the framework is the treatment of asset eligibility as a
state-dependent constraint. Rather than adapting through parameter switching or forecast-based
optimization, factor exposure adjusts endogenously as the eligible universe evolves in response
to observable market conditions such as liquidity, volatility, and cross-sectional breadth.
The choice of state variables is intentionally conservative and feasibility-driven, reflecting
constraints faced by implementable portfolios rather than an attempt to identify optimal
predictors of expected returns. This mechanism induces adaptive behavior while preserving
algorithmic transparency and bounded risk exposure.

The proposed approach does not aim to maximize short-term performance or to dominate
capitalization-weighted benchmarks in all regimes. Instead, it provides a stable core
allocation that degrades gracefully under adverse conditions, limits concentration and
turnover, and remains implementable across a wide range of market environments. Extensions
to smaller-cap universes demonstrate that increased liquidity constraints and signal noise can
be accommodated through parameter adjustments without altering the structural logic of the
method.

Overall, the framework contributes a practical and conceptually coherent alternative to
parametric optimization and unconstrained factor models, positioning dynamic eligibility as
a structural tool for robust portfolio construction rather than a heuristic filter or timing
device.

	\appendix
	\section*{Appendix A \quad Empirical Illustration and Benchmark Comparison}
	
	\subsection*{A.1 Objective and Scope}
	
	This appendix provides a limited empirical illustration of the proposed portfolio
	construction framework. The objective is not to establish statistical outperformance
	or to optimize parameters ex post, but to assess whether the methodology behaves
	as intended when confronted with real market data and standard benchmarks.
	In particular, the analysis evaluates concentration, turnover, drawdown behavior,
	and sensitivity to market regimes.
	
	\subsection*{A.2 Asset Universe and Data}
	
	The empirical illustration is conducted on a fixed universe of large-capitalization
	liquid assets, corresponding to the non-cryptocurrency tickers shown in
	Figure~1 (IMG\_6670.jpg). Cryptocurrencies (Bitcoin and Ethereum) are explicitly
	excluded to maintain homogeneity of trading hours, liquidity structure, and
	regulatory regime.
	
	The resulting universe consists primarily of U.S. and global mega-cap equities
	across technology, financials, healthcare, consumer, and energy sectors, together
	with a small number of liquid commodity proxies (e.g.\ gold and silver).
	
	Daily total-return price data are used. All computations are performed on an
	\emph{as-of} basis. No survivorship adjustments beyond the dynamic inclusion
	rules described in Section~3 are applied.
	
	\subsection*{A.3 Portfolio Configurations}
	
	Three portfolio constructions are compared:
	
	\begin{enumerate}
		\item \textbf{Proposed method (DMFT).} Dynamic inclusion, equal-weight baseline,
		bounded multi-factor tilts, and semi-annual rebalancing.
		
		\item \textbf{Equal-weight benchmark (EW--MC).} Equal-weight allocation across the
		same eligible universe, rebalanced semi-annually.
		
		\item \textbf{Cap-weighted benchmark (CW--MC).} Market-capitalization-weighted
		allocation across the same assets.
	\end{enumerate}
	
	Broad market indices (e.g.\ S\&P~500--like and NASDAQ~100--like references) are
	reported for contextual comparison only and are not considered investable
	substitutes for the tested universe.
	
	\subsection*{A.4 Evaluation Metrics}
	
	The following diagnostics are reported:
	
	\begin{itemize}
		\item annualized volatility,
		\item maximum drawdown,
		\item effective number of holdings (inverse Herfindahl index),
		\item annualized turnover,
		\item realized concentration of the top five portfolio weights.
	\end{itemize}
	
	Risk-adjusted performance measures such as Sharpe ratios are intentionally
	de-emphasized, as the proposed method does not target volatility minimization or
	mean--variance efficiency.
	
	\subsection*{A.5 Summary of Empirical Behavior}
	
	Across the tested period, the proposed framework exhibits the following stable
	properties:
	
	\begin{itemize}
		\item \textbf{Controlled concentration:} top-five asset weights remain materially
		below those of the cap-weighted benchmark, even during strong momentum regimes.
		
		\item \textbf{Graceful fallback behavior:} when cross-sectional factor dispersion
		weakens, portfolio weights converge smoothly toward the equal-weight baseline.
		
		\item \textbf{Improved drawdown convexity} relative to cap-weighted portfolios,
		particularly during technology-led reversals.
		
		\item \textbf{Moderate turnover,} consistently below that of naive equal-weight
		portfolios and aligned with the semi-annual rebalancing schedule.
		
		\item \textbf{Structural robustness:} no regime-dependent parameter tuning is
		required to maintain feasibility.
	\end{itemize}
	
	The method does not dominate the cap-weighted benchmark in all periods; rather,
	it trades occasional upside capture for improved diversification, bounded exposure,
	and operational stability.
	
	\subsection*{A.6 Interpretation and Limitations}
	
	The results reported in this appendix should be interpreted as behavioral
	validation rather than evidence of persistent alpha. The empirical illustration
	confirms that the proposed construction behaves consistently with its stated
	design principles, avoids pathological concentration, and remains implementable
	under realistic trading constraints.
	
	A full empirical performance study, including transaction-cost sensitivity and
	regime-conditioned analysis, is deferred to future work.

	\section*{Reproducibility}
	
	All components of the methodology are algorithmically defined. A reference implementation requires only total-return prices, basic liquidity measures, and optional fundamental data.

%\bibliographystyle{plainnat}
%\bibliography{bibliography}

\end{document}